\documentclass[11pt]{amsart}
\usepackage{amsfonts,amscd}
\usepackage{psfig}
\def\Z{{\mathbb Z}}

\def\P{{\mathbb P}}

\newtheorem{theorem}{Theorem}[section]

\newtheorem{proposition}[theorem]{Proposition}
\theoremstyle{definition}

\theoremstyle{remark}

\numberwithin{equation}{section}
\begin{document}
\title[Symplectic tori in homotopy $E(1)$'s]
{Symplectic tori in homotopy ${\bf E(1)}$'s}
\author{Stefano Vidussi}
\address{Department of Mathematics, Kansas State University, Manhattan, KS-66506}
\email{vidussi@math.ksu.edu}
\maketitle
\baselineskip 18pt
\noindent {\bf Abstract.} This short note presents a simple construction of
nonisotopic symplectic tori representing the same primitive homology
class in the symplectic $4$-manifold $E(1)_{K}$, obtained by knot
surgery on the rational elliptic surface $E(1) =
\P^{2} \# 9 {\bar \P^{2}}$ with the left-handed trefoil knot $K$.
$E(1)_{K}$ has the simplest homotopy type among simply-connected
symplectic $4$-manifolds known to exhibit such property.
\section{Introduction and statement of the results}
Some interest has been devoted in the recent years to the construction
of nonisotopic symplectic surfaces representing the same homology class $\alpha$
of a symplectic $4$-manifold $M$. After the ground-breaking paper of
Fintushel and Stern, that illustrated the existence of this phenomenon for tori
(see \cite{FS2}) and gave a general recipe for the construction of examples,
several authors have variously completed and improved the results
of that paper, see \cite{EP1}, \cite{EP2}, \cite{FS3}, \cite{V2}, \cite{V3}.
As with many other problems of $4$-dimensional topology, one of the open questions
on this subject is how simple can we make a pair $(M,\alpha)$ that presents
such phenomenon. For comparison, it is expected (and partially confirmed in
\cite{Ti}) that for $(\P^{2},\alpha)$ such situation cannot happen.

The complexity of $M$ can be measured in terms of the geography of symplectic
$4$-manifolds, while the complexity of $\alpha$ in terms of its divisibility.
Improvements in this direction have been the result of minor or major modifications of the
original construction of \cite{FS2}. The paper \cite{FS2} presents examples for
$E(1)_{K}$ (with $K$ any fibered knot), for classes with divisibility $4$ or higher.
The first examples representing a primitive class, but on $E(3)$,
appear in \cite{V2}.
Successive results have improved upon these results: examples for $E(2)_{K}$,
for a primitive class, appear in \cite{V3} for infinitely many nontrivial
fibered knots and for all nontrivial fibered knots in \cite{FS3};
similar results for the complex surface $E(2)$ have been later obtained
in \cite{EP2}. Examples of primitive classes for symplectic
manifolds homotopic to $E(1)$ appeared more elusive, as an analysis of the
techniques of the aforementioned papers shows.

In this note we want to improve, at least in
part, this situation, by showing that for infinitely many nontrivial
fibered knots $K$ a primitive class of $E(1)_{K}$ is represented by infinitely
many nonisotopic symplectic tori. The precise statement is the following:
\begin{theorem} \label{infinite} Let $E(n)_{K}$ be the symplectic $4$-manifold
(homeomorphic to $E(n)$) obtained by knot surgery on the elliptic
surface $E(n)$ with the left-handed trefoil knot $K$;
then there exists a nontrivial primitive homology class $[F]$ that
can be represented by infinitely many mutually nonisotopic symplectic tori.
\end{theorem}

This theorem extends without modification to any fibered knot containing $K$ as
a connected summand, as well as to any knot surgery manifold obtained by using the
knot $K$ and several others. Similarly, without much effort, the result can be extended to cover
any multiple of the class $[F]$. Unfortunately, nothing is known to the author
for the case of $E(1)$.

In order to obtain the result of Theorem  \ref{infinite} we will need to
develop a new construction of symplectic tori, which is strongly influenced
by the ideas contained in \cite{FS3} and \cite{V3}.
We will refer to several results and ideas contained in those papers.

({\bf Added in proof:} T. Etg\"u and B.D. Park have announced, in \cite{EP3},
a generalization of Theorem \ref{infinite} to cover the case of any nontrivial
fibered knot and, in \cite{EP4}, a construction of nonisotopic symplectic tori
for the fiber class of $E(1)$. The latter result presents particular interest,
in light of the content of \cite{Ti}.)

\section{Inessential lagrangian tori and essential symplectic tori}
\begin{figure}[h]
\centerline{\psfig{figure=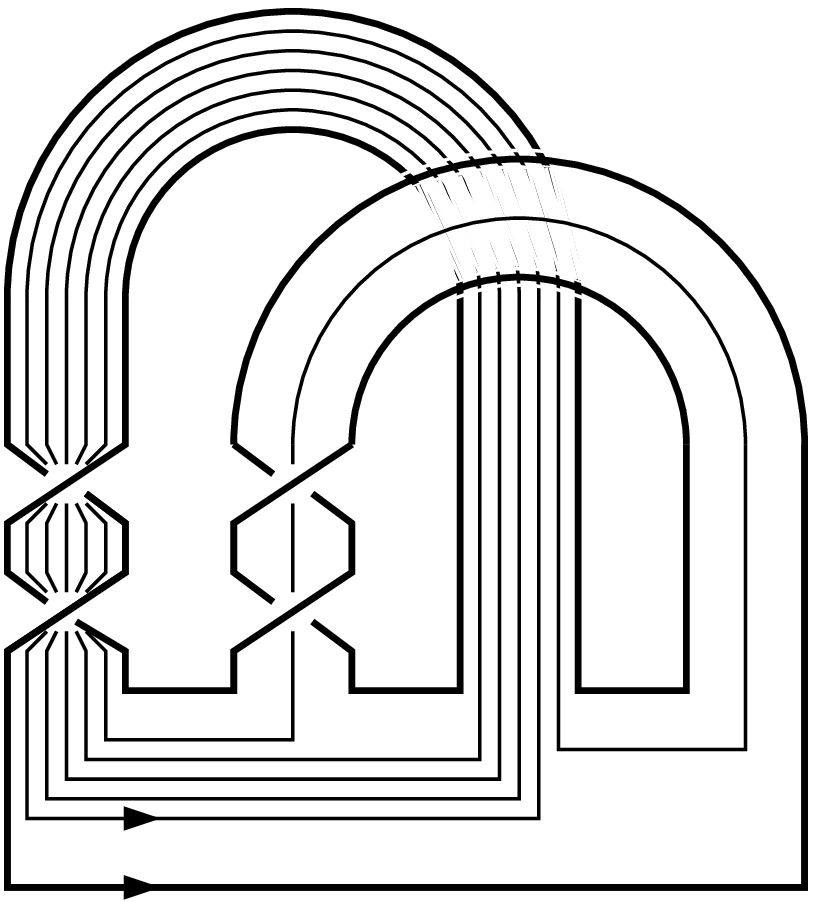,height=38mm,width=38mm,angle=0}
\hspace*{2cm} \psfig{figure=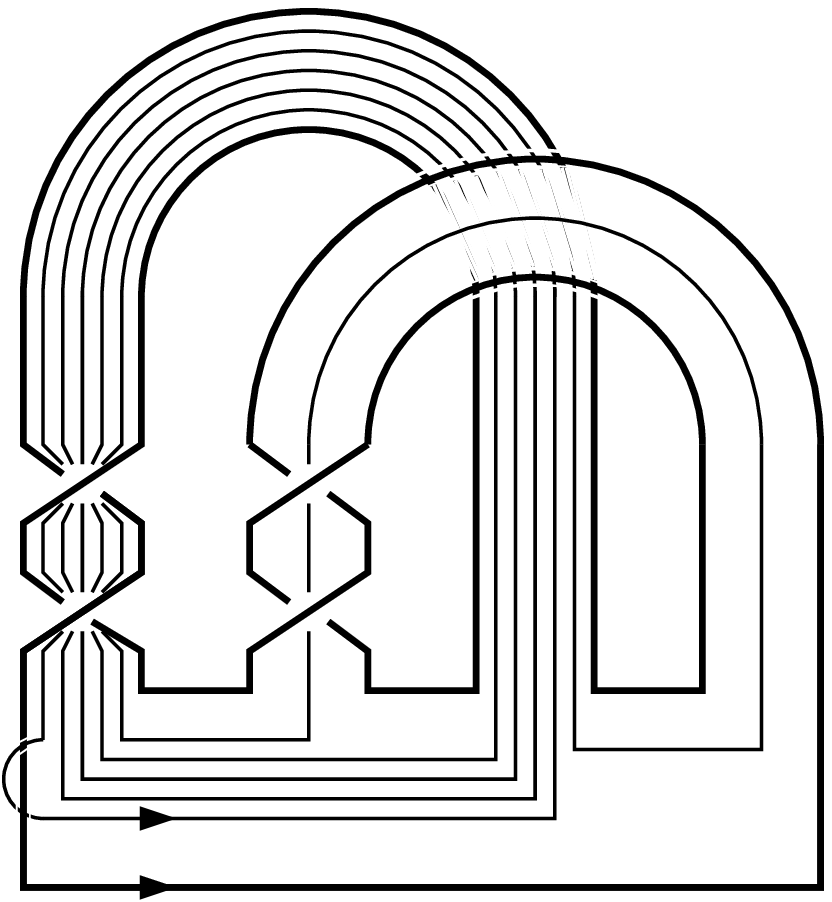,height=38mm,width=40mm,angle=0}}
\caption{\label{trefoil} {\sl The fiber $\Sigma_{K}$ with the curves $\gamma_{p}$
(left) and $\Gamma_{p}$ (right) with $p=5$. }}
\end{figure}

Let $K$ be the left-handed trefoil knot, and $\Sigma_{K}$ its minimal genus
spanning surface, a fixed fiber of the fibration of $S^{3} \setminus \nu K$.
In \cite{V3} the author has constructed a family of nullhomologous lagrangian tori
in the symplectic knot surgery manifold (see \cite{FS1} for definition and
properties)
\begin{equation} E(n)_{K} = (E(n) \setminus \nu F) \cup_{F=S^{1} \times \mu(K)}
S^{1} \times (S^{3} \setminus \nu K) = E(n) \#_{F=S^{1} \times C} S^{1} \times N
\end{equation}  where $N$ is the manifold obtained by $0$-surgery of $S^{3}$
along $K$ and $C$ is the core of the solid torus of the Dehn filling.
The lagrangian tori have the form $S^{1} \times \gamma_{p}$, where
$\gamma_{p}$ is a family of simple closed curves lying on a fiber
$\Sigma_{K} \subset S^{3} \setminus \nu K$, as shown in Figure \ref{trefoil},
l.h.s.
%%%%%%%%%%%%%%
%%%%%%%%%%%%%%

The tori of this family have been proven to be mutually nonisotopic
(see \cite{FS3}). Unfortunately, as these tori are nullhomologous, we cannot make
them symplectic by perturbing the symplectic structure of $E(n)_{K}$.
However, we will obtain a family of symplectic tori by pasting them, in a suitable
sense, with a symplectic torus. We will illustrate this construction
in the following Proposition.
\begin{proposition} There exist a family of closed curves $\Gamma_{p} \in S^{3} \setminus
\nu K$ transverse to the fibration, homologous
to the meridian of $K$ and mutually nonisotopic.
The images $F_{p} \subset E(n)_{K}$ of the tori
$S^{1} \times \Gamma_{p} \subset S^{1} \times
(S^{3} \setminus \nu K)$ define a family of symplectic, framed tori
with self-intersection $0$, representing the class $[F]$. \end{proposition}
\begin{proof} Start by parameterizing smoothly the curve $\gamma_{p}$ as
$\gamma_{p}(s)$, $s \in [-1,1]$, coherently with the orientation.
Consider an annulus $\gamma_{p} \times [-1,1]$ transverse
to the fibration, and hitting the fixed fiber $\Sigma_{K}$ in
$\gamma_{p} = (\gamma_{p},0)$. The fibration restricts, on the annulus, to
a fibration in circles, whose fibers are copies $(\gamma_{p},t)$ of the
knot $\gamma_{p} \subset S^{3}$ pushed in the normal direction.
The annulus $\gamma_{p} \times [-1,1]$ is the image of the diffeomorphism
\begin{equation} \begin{array}{c} \gamma_{p} \times \mbox{Id}: \frac{[-1,1]}{-1=1}
\times [-1,1] \rightarrow \gamma_{p} \times [-1,1] \subset S^{3} \\ \\
(s,t) \mapsto (\gamma_{p}(s),t)
\end{array}  \end{equation}
mapping the standard annulus to a knotted annulus in $S^{3}$.
Define now the curve ${\hat \gamma}_{p}$ as the image
of the ``diagonal", i.e. the curve $(\gamma_{p}(s),s)$. Note that this curve is not
closed and that it is transverse to the fibration of $S^{3} \setminus \nu K$.
Let $M$ be a copy of the meridian to $K$. Up to isotopy, we can
assume that $M$ is transverse to the fibration and that it intersects the annulus
$\gamma_{p} \times [-1,1]$ in the interval $\gamma_{p}(\pm 1) \times [-1,1]$.
We can now define a closed curve $\Gamma_{p}$ by gluing
${\hat \gamma}_{p}$ to the curve $M \setminus (M \cap \gamma_{p} \times [-1,1])$
and smoothing suitably the corners, to make it a smooth curve transverse to
the fibration. (The reader will notice the similarities of this construction
with the circle sum defined in \cite{FS3}, although the purpose here is somewhat
different.) We endow $\Gamma_{p}$ of the canonical framing.
Up to isotopy, the resulting curve (in $S^{3} \setminus \nu K)$) is
drawn in Figure \ref{trefoil}, r.h.s. (Due to difficulties with the drawing, we
have not attempted to reproduce the transversality).
Note that this construction can be carried out simultaneously for all values
of $p$.
As the linking number of $\Gamma_{p}$ and $K$ is equal to $1$ (the linking
number of $\gamma_{p}$ and $K$ is zero as $\gamma_{p} \subset \Sigma_{K}$)
all the curves $\Gamma_{p}$ are homologous to $\mu(K)$ in $H_{1}(S^{3} \setminus
\nu K,\Z)$. However, at least for $p > 1$, they are not isotopic to it;
in fact, as knots in $S^{3}$, $\Gamma_{p}$ is isotopic to $\gamma_{p}$, and
the latter is the torus knot $T_{p,p+1}$, as proven in \cite{V3}.

The torus $S^{1} \times \Gamma_{p}$, that inherits a natural framing coming
from the framing of $\Gamma_{p}$ in $S^{3}$,
is contained in the symplectic manifold $S^{1} \times N$ endowed with the
symplectic form
$dt \wedge \alpha + \epsilon \beta$ (where $\alpha$ is a $1$-form on $N$ defining
the fibration and $\beta$ is a $2$-form on $N$ that restricts to the volume form
on the fibers). As $\Gamma_{p}$ is transverse to the fibers,
$S^{1} \times \Gamma_{p}$ is symplectic, and so is its image in $E(n)_{K}$.
The statement about the homology class follows from the fact that
$[\Gamma_{p}] = [\mu(K)]$ and the identification of $S^{1} \times \mu(K)$
with $F$. \end{proof}

Observe that $[F]$ is a primitive class in $E(n)_{K}$, as $F$ has one intersection
point with the (symplectic) surface obtained by capping  off $\Sigma_{K}$ with a
disk section of $E(n) \setminus \nu F$.

At this point we are ready to prove that infinitely many tori of the family
$F_{p}$ are nonisotopic. Our result is contained in the following theorem, that
implies Theorem \ref{infinite}.
\begin{theorem} Infinitely many elements of the family of homologous symplectic tori
$\{F_{p}: p \geq 1\}$ are not smoothly isotopic. \end{theorem}
\begin{proof} In order to prove this, we observe that an isotopy between two tori
$F_{a}$ and $F_{b}$ implies the existence of a diffeomorphism of the pairs
$(E(n)_{K},F_{a})$ and $(E(n)_{K},F_{b})$. If such a diffeomorphism exists,
the manifolds obtained by summing a copy of $E(1)$ to $E(n)_{K}$ along the two
tori must be diffeomorphic (note that the choice of the gluing map is irrelevant,
because of the use of $E(1)$, see \cite{GS}). The statement follows if
we are able to prove that infinitely many manifolds of the family
\begin{equation} X_{p} = E(n)_{K} \#_{F_{p}=F_{E(1)}} E(1) \end{equation} are not
diffeomorphic. To prove this, we will use a (by now) standard argument on the
number of basic classes (see \cite{V2}). First, note that the manifold
$X_{p}$ is diffeomorphic to a generalized link surgery manifold (see \cite{FS1} and
\cite{V1}), obtained from the $2$-component link $L_{p} = K \cup \Gamma_{p}$:
\begin{equation} X_{p} = E(n,1;L_{p}) := (E(n) \setminus F_{E(n)}) \cup_{T^{3}} S^{1} \times
(S^{3} \setminus \nu L_{p}) \cup_{T^{3}} (E(1) \setminus F_{E(1)}) \end{equation}
where the first gluing map identifies $F_{E(n)}$ with
$S^{1} \times \mu(K)$ and a meridian to $F_{E(n)}$ with $- \lambda(K)$, and the
second gluing map identifies $F_{E(1)}$ with $S^{1} \times \lambda(\Gamma_{p})$ and a meridian to
$F_{E(1)}$ with $\mu(\Gamma_{p})$. The Seiberg-Witten invariant of this manifold
is computed, in terms of the symmetrized Alexander polynomial
$\Delta_{L_{p}}(x,y)$ of $L_{p}$, in \cite{FS1} and \cite{T}:
\begin{equation} SW(X_{p}) = (t_{K} - t^{-1}_{K})^{n-1}
\Delta_{L_{p}}(t^{2}_{K},t^{2}_{\Gamma_{p}}) \end{equation} where $t_{K}$ and $t_{\Gamma_{p}}$
are the elements of $H^{2}(X_{p},\Z)$ Poincar\'e dual to the images of
$S^{1} \times \mu(K)$ and $S^{1} \times \mu(\Gamma_{p})$ respectively (note that
the latter class in nontrivial in $X_{p}$).
As observed in \cite{V2} (and obviously when $n=1$), we can find a lower bound
to the number of basic classes of
$X_{p}$ in terms of the number of nonzero terms in any specialization of
the Alexander polynomial $\Delta_{L_{p}}(x,y)$.
In particular, using Torres formula, we can explicitly
compute the specialization for $x=1$: \begin{equation} \Delta_{L_{p}}(1,y) = \frac{y^{lk(K,\Gamma_{p})}-1}{y-1}
\Delta_{\Gamma_{p}}(y) = \Delta_{\Gamma_{p}}(y). \end{equation}
The polynomial $\Delta_{\Gamma_{p}}(y)$ can be written down explicitly because,
as pointed out above, $\Gamma_{p}$ is the torus knot $T_{p,p+1}$.
In Lemma 6.3 of \cite{V3} the author has shown that the number of nonzero terms of
that polynomial is bounded from below by $p$. Because of that, the number of basic
classes of the manifolds $X_{p}$ grows with no bound with $p$, and therefore
infinitely many of those manifolds are not diffeomorphic. \end{proof}

It is interesting to note in the previous examples the role played by the
nullhomologous lagrangian tori. Comparing this with the constructions of
\cite{FS3} and \cite{V3} we see that in some sense the nonisotopic
symplectic (and lagrangian) essential tori exhibited in these papers
differ by nullhomologous lagrangian tori.

\end{document}